\documentclass[11pt]{article}

\usepackage[alphabetic]{amsrefs}
\usepackage{amsmath,amssymb,enumerate,tikz}
\usepackage[all,cmtip]{xy}
\usepackage[applemac]{inputenc}

\def \1{\mathds{1}}
\def \al{\alpha}
\def \be{\beta}

\def \CF{{\cal F}}

\def \CP{{\cal P}}

\def \df{\ \begin{array}{c} _{\rm def}\\ ^{\displaystyle =}\end{array}\ }
\def \ds{\displaystyle}
\def \e{\emph}
\def \eps{\varepsilon}
\def \Ga{\Gamma}
\def \ga{\gamma}

\def \Hom{\mathrm{Hom}}
\def \Id{{\rm Id}}
\def \Im{\operatorname{Im}}

\def \mathqed{\tag*{$\square$}}
\def \N{{\mathbb N}}
\def \ol{\overline}
\def \om{\omega}

\def \Per{\mathrm{Per}}
\def \ph{\varphi}
\def \prf{{\bf Proof: }}
\def \Q{{\mathbb Q}}
\def \qed{\hfill $\square$}
\def \R{{\mathbb R}}

\def \Sh{\mathrm{Sh}}

\def \Z{{\mathbb Z}}
\def \({\left(}
\def \){\right)}
\def \={\ =\ }

\newcommand{\stack}
[2]{\genfrac{}{}{0pt}{1}{#1}{#2}}

\newcommand{\tto}
[1]{\stackrel{#1}{\longrightarrow}}

\newtheorem{lemma}{Lemma}[section]
\newtheorem{proposition}[lemma]{Proposition}

\newtheorem{conjecture}[lemma]{Conjecture}
\newtheorem{exmples}[lemma]{Examples}
\newenvironment{examples}{\begin{exmples}\nopagebreak\begin{itemize}\nopagebreak\rm}{\end{itemize}\end{exmples}}
\newtheorem{exmple}[lemma]{Example}

\newtheorem{defi}[lemma]{Definition}
\newenvironment{definition}[0]{\begin{defi}\rm}
{\end{defi}}
\newtheorem{theorem}[lemma]{Theorem}

\begin{document}

\pagestyle{myheadings} \markright{HIGHER LOOP SPACES}

\title{Iterated integrals over higher dimensional loops}
\author{Anton Deitmar \& Ivan Horozov}

\date{}
\maketitle

{\bf Abstract:}
We give a definition of higher dimensional iterated integrals based on integration over membranes.
We prove basic properties of this definition and formulate a conjecture which extends Chen's de Rham Theorem for iterated integrals to the membrane case.

$$ $$

\tableofcontents

\newpage
\section*{Introduction}

The classical, or one-dimensional iterated integral is defined for a smooth manifold $X$ as well as 1-forms $\om_1,\dots,\om_m$ and a path $\ga:I\to X$ by
$$
\int_\ga\om_1\dots\om_m=\int_{t_1<t_2<\dots<t_m}\ga^*\om_1\wedge\dots\wedge\ga^*\om_m.
$$
See \cite{Ch} for more on this.
One of the most important results of Chen's is the so called de Rham Theorem of Chen, which we state now.
Fix a base point $x_0\in X$ and let $B_s(X)_{x_0}$ denote the real vector space of all maps from the set of loops at $x_0$ to $\R$ which are linear combinations of  iterated integrals of length at most $s$.
Let $B_s(X)_{x_0}^{hom}$ denote the subspace of homotopy-invariant maps.
Chen's de Rham Theorem says that the evaluation map
$$
B_s(X)_{x_0}^{hom}\to \Hom_\Z(\Z[\pi_1(X,x_0)]/J^{s+1},\R)
$$
is an isomorphism of real vector spaces, where $\Z[\pi_1(X,x_0)]$ is the group ring and $J$ the augmentation ideal.

In this paper we propose a higher dimensional version of iterated integrals.
We investigate their basic properties and state the conjecture that the resulting map
$$
B_s^n(X)_{x_0}^{hom}\to \Hom_\Z(\Z[\pi_n/\pi_n^{<n}]/J^{s+1},\R)
$$
is an isomorphism.
here $\pi_n$ is the higher homotopy group and $\pi_n^{<n}$ is the subgroup generated by maps from spheres of lower dimension.
We prove this conjecture in the case $s=1$.

Similar, but less general iterated integrals over higher dimensional membranes were first defined in \cite{H}, where they are applied to Hilbert modular forms. They also are used in a construction of Multiple Dedekind zeta functions \cite{Multi}. The definition of such iterated integrals was done in a search of a higher dimensional analogue of Manin's non-commutative modular symbol \cite{M}.

\section{Iterated integrals on a membrane}

Let  $n\in\N$.
Instead of paths we consider membranes.
A \e{$n$-dimensional membrane} on $X$ is a continuous map $g:I^n\to X$, which is continuously differentiable outside a closed set $S\subset I^n$ of measure zero.
On the set $I^n$ we introduce the partial order
$$
t<s\quad\Leftrightarrow\quad t_1<s_1,\dots,t_n<s_n.
$$
The direct generalization of the one-dimensional iterated integral would be this:
For $n$-forms $\om_1,\dots,\om_r$ we define the \e{$n$-dimensional iterated integral} as
$$
\int_g\om_1\dots\om_r=\int_{t_1<\dots<t_r}g^*\om_1\wedge\dots\wedge g^*\om_r.
$$
In this paper we present a more general approach.
The following definition of an $n$-dimensional 
iterated integral is motivated by physics. 
Consider a number of
\begin{itemize}
\item {\bf $s$ events}, watched by
\item {\bf $n$ observers}
\end{itemize}
 in special relativity.  Let $t_\nu$ be 
the time coordinate for observer $\nu$, where $\nu=1,\dots,n$. Now consider $s$ events, $e_1,\dots,e_s$. The 
simplest case is, when each observer sees the events in the same order. However, time is relative, so it is 
possible that each observer sees the events in a different order. Let $\rho_\nu$ be a permutation of $s$ 
elements, describing the order in which observer $\nu$ sees the events. 
For each $\nu=1,\dots,n$, we cut the admissible interval for the time variable $t_\nu$ into subintervals, where the different 
events occur.
For fixed $\nu$, let 
$$0<t_\nu^1<t_\nu^2<\dots<t_\nu^m<1$$
values of $t_\nu$ in each subinterval.

Consider variables by $t^\sigma_\nu$ for $\nu=1,\dots, n$ and $\sigma=1,\dots,s$. 
The subscript $\nu$ corresponds to $\nu$-th 
time direction for the observer $\nu$. 
 
Let $D$ be a domain in terms of the  variables $t^\sigma_\nu$, defined by
$$D=\{(t_\nu^\sigma)|0<t_\nu^1<t_\nu^2<\dots<t_\nu^s<1,\ \nu=1,\dots,n\}\subset I^{sn},$$
where $I=[0,1]$ is the unit interval.

We associate differential $n$-forms $\omega_\sigma$ on a manifold $X$ to each event $\sigma$, for
$\sigma=1,2,\dots,s$.
Let 
$$g:[0,1]^n\rightarrow X$$
be a membrane. 
The observer $\nu$ sees the event  $\sigma$, realized as a differential form $g^*\omega_\sigma$ at time $t_\nu^{\rho_\nu(\sigma)}.$ 
Therefore $g^*\omega_\sigma$ depends on 
$(t_1^{\rho_1(\sigma)},\dots,t_n^{\rho_n(\sigma)})$, which  are the time coordinates for each observer, in which they 
see the event $\sigma$.

To each $\sigma=1,\dots,s$ and the permutations $\rho_1,\dots,\rho_n$ we associate a map
\begin{align*}
\phi_{\sigma,\rho}:I^{sn}&\to I^n\\
(t_\nu^\sigma)_{\sigma,\nu}&\mapsto \(t_1^{\rho_1(\sigma)},\dots,t_n^{\rho_n(\sigma)}\).
\end{align*}

\begin{definition}
\label{hIteratedIntegral}
An $n$-dimensional iterated integral in terms of $n$-forms $\omega_\sigma$, $\sigma=1,\dots,s$, 
permutations $\rho_\nu$, for $\nu=1,\dots,n$ and a smooth map $g$, is defined as
$$
\int_g^\rho \om_1\cdots\om_s \df
\int_D\phi_{\rho,1}^*\om_1\wedge\dots\wedge\phi_{\rho,s}^*\om_s.
$$
\end{definition}

In the case when $\rho_\nu=\Id$ for all $\nu$, we leave out the superscript $\rho$.

\begin{examples}
\item In the classical case $n=1$ of one observer there is only one permutation $\rho\in\Per(s)$ and one gets
$$
\int_g^\rho\om_1\cdots\om_s=\det(\rho)\int_g\om_{\rho(1)}\cdots\om_{\rho(s)}.
$$
\item In the special case of $2$ observers one may use the following type of diagram for visualization.
We have the observers $A$ and $B$, and we draw a time axis for each of them. Say there are three events $e_1,e_2,e_3$, then the diagram
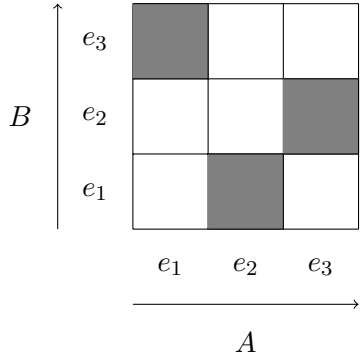
\begin{center}
\begin{tikzpicture}
\draw[step=1cm] (0,0) grid (3,3);
\draw[fill=black!50](1,0)--(2,0)--(2,1)--(1,1);
\draw[fill=black!50](2,1)--(3,1)--(3,2)--(2,2);
\draw[fill=black!50](0,2)--(1,2)--(1,3)--(0,3);
\draw (0.5,-.5)node{$e_1$};
\draw (1.5,-.5)node{$e_2$};
\draw (2.5,-.5)node{$e_3$};
\draw (-0.5,0.5)node{$e_1$};
\draw (-.5,1.5)node{$e_2$};
\draw (-.5,2.5)node{$e_3$};
\draw[->](0,-1)--(3,-1);
\draw[->](-1,0)--(-1,3);
\draw(1.5,-1.5)node{$A$};
\draw(-1.5,1.5)node{$B$};
\end{tikzpicture}
\end{center}
means that observer $A$ sees the events in the order $e_3,e_1,e_2$ and $B$ sees them as $e_2,e_3,e_1$.
The shaded region indicates, where the integration takes place.
Note that these diagrams only cover the case $\rho_1\rho_2=\Id$.
\end{examples}

There is also a different way of presenting the iterated integral as follows.
Define
$\phi_\rho:I^{sn}\to I^{sn}$ by
$$
\phi_\rho(t)_\nu^\sigma=t_\nu^{\rho(\sigma)}.
$$
Let $p_\sigma:X^s\to X$ be the $\sigma$th projection and consider the $sn$-form on $X^s$,
$$
\om=p_1^*\om_1\wedge\dots\wedge p_s^*\om_s.
$$
We denote the map $(I^n)^s\to X^s$ with coordinates $g$ also by $g$.

\begin{lemma}
We have
$$
\int_g^\rho\om_1\cdots\om_s=\int_D(g\circ\phi_\rho)^*\om=\int_{D_\rho}g^*\om=\int_{g(D_\rho)}\om.
$$
here $D_\rho$ is the set of all $t\in I^{sn}$ with $0<t_\nu^{\rho_\nu(\sigma)}<t_\nu^{\rho_\nu(\sigma+1)}<1$ for all $\sigma,\nu$.
\end{lemma}

\prf Clear.
\qed

\section{Reparametrization}
Let $\ph:I^n\to I^n$ be a piecewise diffeomorphism, that is $\ph$ is a homeomorphism, such that $\ph$ and $\ph^{-1}$ are continuously differentiable outside a closed set of measure zero.
We say that $\ph$ is \e{monotonic}, if
$$
x\le y\quad\Rightarrow\quad \ph(x)\le \phi(y),
$$
where we say $x\le y$ if $x_j\le y_j$ for every $1\le j\le n$.
We define $\ph^s:(I^n)^s\to (I^n)^s$ by
$$
\ph^s(t_1,\dots,t_s)=(\ph(t_1),\dots,\ph(t_s)).
$$

\begin{lemma}
If $\ph$ is a monotonic homeomorphism, then 
$$
\ph^s(D_\rho)=D_\rho.
$$ 
\end{lemma}

\prf In our notation, $D$ is the set of all $x\in (I^n)^s$ with $x_1<x_2<\dots<x_s$.
For $x\in D$ we have $\ph(x_1)\le\ph(x_2)\le\dots\le\ph(x_n)$, hence $\ph^s(x)\in \ol D$, the closure of $D$, i.e., $\ph(D)\subset \ol D$.
But as $\ph$ is a homeomorphism, the image $\ph(D)$ is open. Every open subset of $\ol D$ lies in $D$, so $\ph(D)\subset D$.
The inverse map $\ph^{-1}$ is monotonic as well, hence $\ph(D)=D$.
We have $\ph^s\circ\phi_\rho=\phi_\rho\circ\ph^s$, so that finally
\begin{align*}
\ph^s(\phi_\rho(D))&=\phi_\rho(\ph^s(D))=\phi_\rho(D).\mathqed
\end{align*}

\begin{proposition}
\begin{enumerate}[\rm (a)]
\item We have
$$
\int_g^\rho\om_1\cdots\om_s=
\int_{g\circ\ph}^\rho\om_1\cdots\om_s
$$
for every monotonic piecewise diffeomorphism $\ph:I^n\to I^n$.
\item We have
$$
\int_{F\circ g}^\rho \om_1\dots\om_s =
\int_{g}^\rho (F^*\om_1)\dots(F^*\om_s) 
$$
fr every smooth map $F:X\to Y$.
\end{enumerate}\end{proposition}

\prf
(a) By $g_\rho=g^s\circ\phi_\rho$ we get
\begin{align*}
\int_{g\circ\ph}^\rho\om_1\cdots\om_s
&= \int_D (g^s\circ\ph\circ\phi_\rho)^*\om= \int_D (\ph\circ\phi_\rho)^*(g^s)^*\om\\
&= \int_{\ph\circ\phi_\rho(D)} (g^s)^*\om= \int_{\phi_\rho(D)} (g^s)^*\om\\
&= \int_{D} \phi_\rho^*(g^s)^*\om= \int_{D} g_\rho^*\om= \int_{g}^\rho\om_1\cdots\om_s.
\end{align*}

For (b) note $(F\circ g)_\rho=F^s\circ g_\rho$, so that
\begin{align*}
\int_{F\circ g}^\rho \om_1\dots\om_s&= \int_D(F^s\circ g_\rho)^*\om=\int_Dg_\rho^*(F^s)^*\om\\
&=
\int_{g}^\rho (F^*\om_1)\dots(F^*\om_s).
\mathqed
\end{align*}

\section{Homotopy invariance}

In this section we are going to show that for holomorphic differential forms, the iterated integral is invariant under homotopic deformations of $g$, which preserve certain foliations, see below.

For each $\nu=1,\dots,n$ the fibers of the projection $p_\nu:I^n\to I$; $(x_1,\dots,x_n)\mapsto x_\nu$ form the leaves of a foliation $\CP_\nu$.

Let $X$ be a complex manifold of complex dimension $n$ and let $\CF_1,\dots,\CF_n$ be non-trivial foliations with complex submanifolds as leaves.
We say that a membrane $g:T^n\to X$ is \e{$\CF$-admissible}, if $g$ maps $\CP_\nu$-leaves to $\CF_\nu$-leaves.
A homotopy $h:I^n\times I\to X$ is called \e{$\CF$-admissible}, if each intermediate $g_u=h(.,u)$ is $\CF$-admissible.

\begin{theorem}
Let $X$ be a complex manifold of complex dimension $n$.
Let $\CF_1,\dots,\CF_n$ be complex foliations on $X$.
Let $\om_1,\dots,\om_s$ be holomorphic $n$-forms on $X$.

Assume that two $\CF$-admissible membranes $g_0,g_1$ are homotopic with a $\CF$-admissible homotopy.
Then we have
$$
\int_{g_0}^\rho\om_1\cdots\om_s=\int_{g_1}^\rho\om_1\cdots\om_s.
$$
\end{theorem}

\prf
Write $F_i=\int_{g_i}^\rho\om_1\cdots\om_s$.
Let $I$ denote the unit interval $[0,1]$.
Let $\om=p_1^*\om_1\wedge\dots\wedge p_s^*\om_s$.
Then $\om$ is a holomorphic $sn$-form on $X^s$, so $d\om =0$.
Define
$$
h_\rho:D\times I\to X^s
$$ by
$$
h_\rho(t,u)=\(h(\phi_{\rho,1}(t),u),\dots,h(\phi_{\rho,s}(t),u)\).
$$
Note that $h_\rho(t,i)=g_{i,\rho}$ for $i=0,1$.
By Stokes's Theorem we get
\begin{align*}
0&=\int_{D\times I}d(h_\rho^*\om) =\int_{\partial(D\times I)}h_\rho^*\om\\
&= \underbrace{\int_{D\times\{ 1\}}h_\rho^*\om}_{=F_1}-\underbrace{\int_{D\times\{ 0\}}h_\rho^*\om}_{=F_0} +\sum_{\nu=1}^n\sum_{\sigma=0}^s\eps_{\sigma,\nu}\int_{D_{\sigma,\nu}\times I}h_\rho^*\om,
\end{align*}
where $\eps_{\sigma,\nu}\in\{\pm 1\}$ is a sign and $D_{\sigma,\nu}$ is defined as follows.
Put
$t_\nu^0=0$ and $t_\nu^{s+1}=1$ and define
$$
D_{\sigma,\nu}=\left\{ (t_\al^\be)_{\stack{1\le\al\le n}{1\le\be\le s}}: \begin{array}{c}0<t_\al^1<\dots<t_\al^s<1\text{  for }\al\ne\nu,\\
t^0_\nu<\dots<t_\nu^{\sigma}=t_\nu^{\sigma+1}<\dots<t^{s+1}_\nu\end{array}\right\}\subset I^{sn}.
$$

We will show that $\int_{D_{\sigma,\nu}\times I}h_\rho^*\om=0$ for all $\sigma,\nu$.
Note that $\int_{D_{\sigma,\rho}\times I}h_\rho^*\om=\int_{h(D_{\sigma,\rho}\times I)}\om$.
Since each $g_u$ is $\CF$-admissible, $g_{u,\rho}$ maps $D_{\sigma,\nu}$ to the set $X_{\sigma,\nu}$ of all $x_1,\dots,x_s\in X^s$ for which $x_{\rho(\sigma)}$ and $x_{\rho(\sigma+1)}$ lie in the same $\CF_\nu$-leaf.
The latter is a complex submanifold of $X^s$ of dimension $sn-1$.
Since this is true for all $u$, it follows that
$$
h(D_{\sigma,\nu}\times I)\subset X_{\sigma,\nu},
$$
 so that we integrate a holomorphic $sn$-form over a submanifold of dimension $\le sn-1$, where all holomorphic $sn$-forms vanish, hence the integral is zero.
\qed

\section{Shuffle relations}
Let $s,s'\in\N$.
An  \e{$(s,s')$-shuffle} is a permutation $\sigma\in\Per(s+s')$ such that
$$
\sigma(1)<\dots<\sigma(s)\quad\text{and}\quad \sigma(s+1)<\dots<\sigma(s+s').
$$
For two permutations $\rho\in\Per(s)$ and $\rho'\in\Per(s')$, a \e{$(\rho,\rho')$-shuffle} is a permutation $\sigma\in \Per(s+s')$ of the form
$$
\sigma=\tau\circ(\rho,\rho'),
$$
where $\tau$ is an $(s,s')$-shuffle.
Let $\Sh(\rho,\rho')$ denote the set of all $(\rho,\rho')$-shuffles.
A given permutation $\sigma\in\Per(s+s')$ lies in $\Sh(\rho,\rho')$ if and only if 
$$
\sigma(\rho(1))<\dots<\sigma(\rho(s))\quad\text{and}\quad \sigma(n+\rho'(1))<\dots<\sigma(s+\rho'(s')).
$$
Let $\rho=(\rho_1,\dots,\rho_n)$ and $\rho'=(\rho_1',\dots,\rho_n')$ be tuples of permutations in $\Per(s)$ and $\Per(s')$ respectively.
Then define
$$
\Sh(\rho,\rho')=\prod_{\nu=1}^n\Sh(\rho_\nu,\rho_\nu').
$$

\begin{proposition}[Shuffle relations]
We have
$$
\(\int_g^\rho\om_1\dots\om_s\)\(\int_g^{\rho'}\om_{s+1}\dots\om_{s+s'}\)=\sum_{\tau\in\Sh(\rho,\rho')}\int_g^\tau\om_1\dots\om_{s+s'}.
$$
\end{proposition}

\prf
We write the product on the left as
\begin{multline*}
\(\begin{array}{cc}{\ds\int}&\om_1\wedge\dots\wedge\om_s\\
{\begin{array}{ccc}t_1^{\rho_1(1)}<&\dots&<t_1^{\rho_1(s)}\\
\vdots&&\vdots\\
t_n^{\rho_n(1)}<&\dots&<t_n^{\rho_n(s)}\end{array}}
\end{array}\)\times\\
\(\begin{array}{cc}{\ds\int}&\om_1\wedge\dots\wedge\om_s\\
{\begin{array}{ccc}s_1^{\rho_1'(1)}<&\dots&<s_1^{\rho_1'(s)}\\
\vdots&&\vdots\\
s_n^{\rho_n'(1)}<&\dots&<s_n^{\rho_n'(s)}\end{array}}
\end{array}\).
\end{multline*}
Consider points $t,x$ in the respective domains of integration, which are chosen generically in the sense that no coordinates of $t$ match any coordinates of $x$.
Then the order of the points $t_1^{\rho_1(1)},\dots,t_1^{\rho_1(m)},x_1^{\rho_1(1)},\dots,x_1^{\rho_1'(s)}$ determines an $(\rho_1,\rho_1')$-shuffle $\tau_1$. Repeat this in the other rows to get a $(\rho,\rho')$-shuffle $\tau$.
All such shuffles appear and the domain of integration is, up to a set of measure zero, the disjoint union over all $\tau$.
This gives the claim.
\qed

\section{Composition}
We say that a membrane $g$ is \e{closed} with \e{base point} $x_0\in X$ if $g$ maps the entire boundary of $I^n$ to the singleton $\{x_0\}$.

Let the membranes $g_1,g_2$ be closed with the same base point $x_0$.
We define a membrane  or $g_1g_2$ by
$$
g_1g_2(t)=\begin{cases}g_1(2t) & 0\le t_\nu\le\frac12\ \forall_\nu\\
g_2(2t-1)& \frac12< t_\nu\le 1\ \forall_\nu\\
x_0& \text{otherwise.}
\end{cases}
$$

\begin{proposition}
Let $g_1,g_2$ be closed membranes with the same base point.
Then
$$
\int_{g_1g_2}\om_1\dots\om_s=\sum_{j=0}^s\(\int_{g_1}\om_1\dots\om_j\)\(\int_{g_2}\om_{j+1}\dots\om_s\).
$$
\end{proposition}

\prf
The left hand side equals
$$
\int_{t_\nu^1<\dots<t_\sigma^s}(g_1g_2)^*\om_1\wedge\dots\wedge (g_1g_2)^*\om_s.
$$
Let $F$ be the set of all $t\in I^n$ such that either all coordinates are $\le\frac12$ or all coordinates are $\ge \frac12$.
If any of the $t^1,\dots t^m$ is outside $F$, then locally the integral is zero as one of the forms $(g_1g_2)^*\om_\sigma$ is zero.
Therefore, the integral can be restricted to $F^s$ and the order conditions force in that if all coordinates of $t^\sigma$ are $\ge\frac12$, then the same holds for $t_{\sigma+1}$.
This gives the claim.
\qed

Let $P^nX$ denote the set of all membranes $g:I^n\to X$ and write $B_s^n(X)$ for the set of all maps $P^nX\to \R$ which are linear combinations of iterated integrals of lenght $\le s$.

We write $L^n_{x_0}X$ or just $L^nX$ for the space of loops at $x_0$, i.e., the set of $g\in P^nX$ with $g(\partial(I^n))=\{ x_0\}$.
The set of restrictions of elements of $B_s^n(X)$ to $L^nX$ is denoted by $B_s^n(X)_{x_0}$.

We define the space $L^nX_{par}$ to be the quotient of $L^nX$ modulo the equivalence relation given by monotonic reparametrization. 
The composition law $(g_1,g_2)\mapsto g_1g_2$ is easily seen to be associative on $L^nX_{par}$.
This being the case, the composition turns the free abelian group $\Z[L^nX_{par}]$ into a ring.

\begin{proposition}\label{prop5.2}
Every $\om\in B_s^n(X)_{x_0}$ factors over $L_{x_0}^nX_{par}$.
We extend it linearly to a map on the ring $\Z[L^nX_{par}]$.
For $\al_1,\dots,\al_r\in L^nX_{par}$ we set $\eta=(\al_1-1)\cdots(\al_r-1)$.
For $n$-forms $\om_1,\dots,\om_s$ we then have
$$
\int_\eta\om_1\cdots\om_s=\begin{cases}\prod_{j=1}^r\int_{\al_j}\om_j &\text{if } s=r\\
0&\text{if }s<r.\end{cases}
$$
\end{proposition}

\prf
The proof is an iterated application of the composition formula
$$
\int_{g_1g_2}\om_1\dots\om_s=\sum_{j=0}^s\(\int_{g_1}\om_1\dots\om_j\)\(\int_{g_2}\om_{j+1}\dots\om_s\).
$$
Let $a_k=(\al_1-1)\cdots(\al_k-1)\al_{k+1}\cdots\al_r$.
By induction one shows
$$
\int_{a_k}\om_1\cdots\om_s=\sum_{0<i_1<\dots<i_k\le i_{k+1}\le\dots\le i_s\le s}\int_{\al_1}\om_1\cdots\om_{i_1}\dots\int_{\al_r}\om_{i_s+1}\cdots\om_s.
$$
If $k=r$, then the sum becomes empty for $s<r$ and has one term only for $r=s$.
\qed

\section{The conjectural de Rham Theorem}

Let $B_s^n(X)_{x_0}^{hom}$ be the set of all elements of $B_s^n(X)_{x_0}$ which are invariant under homotopies which leave the boundary of $I^n$ fixed.
Proposition \ref{prop5.2} implies that each $\om\in B_s^n(X)_{x_0}^{hom}$ induces a map $\Z[\pi_n(X,x_0)]/J^{s+1}\to \R$, where $J$ is the augmentation ideal in the group ring $\Z[\pi_n(X,x_0)]$.
So we get an injection
$$
B_s^n(X)_{x_0}^{hom}\hookrightarrow \Hom_\Z(\Z[\pi_n(X,x_0)],\R).
$$
In the case $n=1$ Chen's de Rham Theorem says that this map is also onto.
For higher $n$, this can't be the case, because if $n$ is bigger than the dimension of $X$, then the space $B_s^n(X)$ is zero, as every $n$-form vanishes.

{\bf Question:}
\begin{center}
{\it What is the image of $B_s^n(X)_{x_0}^{hom}$ inside the space $\Hom_\Z(\Z[\pi_n(X,x_0)],\R)$?}
\end{center}

We formulate a conjecture: Let 
$\pi_n^{<n}(X,x_0)$ denote the subgroup of $\pi_n(X,x_0)$ generated by the images of all maps 
$$
\al_*:\pi_n(S^k)\to\pi_n(X),
$$
where $\al$ ranges over $\pi_k(X,x_0)$ and $1\le k<n$.

\begin{lemma}
For $n=1$ one has $\pi_n^{<n}(X,x_0)=\pi_n(X,x_0)$.
For $n\ge 2$ the group $\pi_n^{<n}(X,x_0)$ is the kernel of the Hurewicz map $h_n:\pi_n(X)\to H_n(X)$.
One has
$$
\pi_n^{<n}(X)\otimes\Q=\bigcup_{\stack{x_0\in Y\subset X}{\dim Y<n}}\Im(\pi_n(Y,x_0)\to\pi_n(X,x_0))\otimes\Q,
$$
where the union is extended over all connected subsets $Y$, which are finite unions of submanifolds of dimension $<n$.
\end{lemma}

\prf
For each non-trivial $\al\in\pi_k(X,x_0)$ with $k<n$ take $\al$ to be the glueing prescription for glueing a new $(k+1)$-cell onto $X$ and let $W$ denote the resulting CW-complex.
Then $W$ satisfies the conditions for Hurewicz's theorem, so the map $\pi_n(W)\to H_n(W)$ is an isomorphism.
By the long exact homology sequence for the pair $(W,X)$ one gets that the natural map $H_n(X)\to H_n(W)$ is an isomorphism.
Let $\phi:\pi_n(X)\to\pi_n(W)$ be the induced map.
The commutative square
$$
\xymatrix{
\pi_n(X)\ar[r]^{\phi}\ar[d]_{h_n}&\pi_n(W)\ar[d]^\cong\\
H_n(X)\ar[r]^\cong&H_n(W)
}
$$
shows that $\ker(h_n)=\ker(\phi)$.
Let $H$ be the subgroup generated by all images of $\al_*$ as in the lemma. Then $H\subset\ker(\phi)$.
For the converse, let $\al\in\ker(\phi)$. Then there is a contracting homotopy $h$ in $W$.
As the image of $h$ is compact, it hits only finitely many new cells.
By modifying $h$, where it leaves $X$, one gets a homotopy that stays in $X$ and moves $\al$ not to zero, but to an element of $H$.
This proves the first assertion.

For the second, note that 
$$
\bigcup_{\stack{x_0\in Y\subset X}{\dim Y<n}}\Im(\pi_n(Y,x_0)\to\pi_n(X,x_0))
$$
is a subgroup of $\pi_n$ which contains $H=\pi_n^{<n}$ and which maps to zero under the Hurewicz map followed by the de Rham map 
\begin{align*}
\pi_n(X)\otimes\Q\to H_n(X)\otimes\Q\hookrightarrow H^n(X,\R)^*.
\mathqed
\end{align*}

\begin{conjecture}[Higher dimensional de Rham Theorem]
We conjecture that iterated integration gives an isomorphism
$$
B_s^n(X)_{x_0}^{hom}\tto\cong \Hom_\Z(\Z[\pi_n/\pi_n^{<n}]/J^{s+1},\R),
$$
where $J$ is the augmentation ideal of the group ring.
\end{conjecture}

For $n=1$, this is Chen's Theorem. For $s=1$ we prove it below.

\begin{proposition}
Let $n\ge 2$ and let $X$ be a smooth and compact manifold.
Then the conjecture holds for $s=1$. That is, we have
$$
B_1^n(X)_{x_0}^{hom}\tto\cong \Hom_\Z(\Z[\pi_n/\pi_n^{<n}]/J^{2},\R).
$$
\end{proposition}

\prf
Observe first, that a single  $n$-form $\om$ is homotopy invariant (either as element of $B_1^n(X)$ or of $B_1^n(X)_{x_0}$), if and only if the form $\om$ is closed.
Since $\Ga=\pi_n/\pi_n^{<n}$ is abelian, the map $\Ga\to J/J^2$; $\ga\to [\ga-1]$ is an isomorphism.
Therefore,
\begin{align*}
\Hom_\Z(\Z[\Ga]/J^2,\R) &=\Hom_\Z(\Z\oplus J/J^2,\R)\\
&=\R\oplus\Hom_\Z(\Ga,\R)\\
&= \R\oplus\Hom_\Z(\Im(h_n),\R).
\end{align*}
The  pairing between homology and de Rham cohomology is given by integration and it identifies $H_{dR}^n(X,\R)$ with $\Hom_\Z(H_n(X),\R)$.
Recall that the space $B_1^n(X)_{x_0}^{hom}$ is the direct sum of $\R$ and  the set of restrictions of elements of $\Hom_\Z(H_n(X),\R)$ to $\Im(h_n)$.
The proposition follows.
\qed

{\small Anton Deitmar, Mathematisches Institut,
Auf der Morgenstelle 10,
72076 T\"ubingen,
Germany,
\tt deitmar@uni-tuebingen.de}

{\small Ivan Horozov,
Department of Mathematics,
Washington University in Saint Louis,
One Brookings Drive, St. Louis, MO 63130
\tt horozov@math.wustl.edu}
\end{document}